# Most Probable Phase Portraits of Stochastic Differential Equations and Its Numerical Simulation

**Bing Yang[1], Zhu Zeng[2] and Ling Wang[3]**

[1] *School of Mathematics and Statistics, Huazhong University of Science and Technology, email: yang_bobby@qq.com*
[2] *School of Mathematics and Statistics, Central China Normal University, e-mail: zhuzeng_style@qq.com*
[3] *School of Mathematics and Statistics, Central China Normal University, e-mail: wanglingccnu@qq.com*

**Abstract.** A practical and accessible introduction to most probable phase portraits is given. The reader is assumed to be familiar with stochastic differential equations and Euler-Maruyama method in numerical simulation. The article first introduce the method to obtain most probable phase portraits and then give its numerical simulation which is based on Euler-Maruyama method. All of these are given by examples and easy to understand.

## 1. Introduction

A phase portrait is a geometric representation of the trajectories of a dynamical system in the phase plane. For deterministic dynamical systems, phase portraits provide geometric pictures of dynamical orbits, at least for lower dimensional systems. However, a stochastic dynamical system is quite different from the deterministic case. There have been some options of phase portraits already. But they are all limited in some ways. Hence, in this article we explain a new kind of phase portraits — most probable phase portraits, which is first proposed by Prof. Duan in his recent published book [1, §5.3].

In the next section, we introduce you the history of phase portraits for SDEs, including some examples to show how to get most probable phase portraits. In section 3 we explain our motivation in this article is to establish the numerical method of most probable phase portraits. In section 4, we show our numerical method through examples and compare it with the real result. Finally, we give a summary of our article.

## 2. History

### 2.1 Earlier methods in phase portraits
There are two apparent options of phase portraits, which are mean phase portraits and almost sure phase portraits.

#### 2.1.1 Mean phase portraits
Let us consider a simple linear SDE system

$$dX_t = 3X_t + dB_t \tag{2.1}$$

The mean $EX_t$ evolves according to the linear deterministic system

$$\frac{d}{dt}EX_t = 3EX_t \tag{2.2}$$

which is the original system without noise. In other words, the mean phase portrait will not capture the impact of noise in this simple linear SDE system.

The situation is even worse for nonlinear SDE system. For example, consider

$$dX_t = (X_t - X_t^3)dt + dB_t \tag{2.3}$$

Take mean on both side of this SDE to get

$$\frac{d}{dt}EX_t = EX_t - E(X_t^3) \tag{2.4}$$

Thus, we do not have a 'closed' differential equation for the evolution of mean $EX_t$, because $E(X_t^3) \neq (EX_t)^3$. This is a theoretical difficult for analyzing mean phase portraits for stochastic system. The same difficulty arises for mean-square phase portraits and higher-moment phase portraits.

#### 2.1.2 Almost sure phase portraits
Another possible option is to plot sample solution orbits for an SDE system, mimicking deterministic phase portraits. If we plot representative sample orbits in the state space, we will see it could hardly offer useful information for understanding dynamics. Moreover, each sample orbit is a possible 'outcome' of a realistic orbit $X_t$ of the system. But which sample orbit is most possible or maximal likely? This is determined by the maximizers of the probability density function $p(x,t)$ of $X_t$, at every time t.

### 2.2 New method in phase portraits -- Most probable phase portraits
In the last section, we discussed two sorts of phase portraits. Mean phase portraits and almost sure phase portraits. Mean phase portraits has difficulties for nonlinear SDE systems and higher-moment phase portraits. Almost sure phase portraits shows us a very complicated picture. It is difficult to find useful information.

In this section, we introduce you a deterministic geometric tool —most probable phase portraits, which is first proposed by Professor J. Duan, see [1, §5.3]. The most probable phase portraits provide geometric pictures of most probable or maximal likely orbits of stochastic dynamical systems. It is based on Fokker-Planck equations.

For an SDE system in $R^n$

$$dX_t = b(X_t)dt + \sigma(X_t)dB_t \quad , \quad X_0 = \xi \tag{2.5}$$

The Fokker-Planck equation for the probability density function $p(x,t)$ of $X_t$ is

$$\frac{\partial}{\partial t} p(x,t) = A^* p(x,t) \tag{2.6}$$

With initial condition $p(x,0) = \delta(x-\xi)$. Recall that the Fokker-Planck operator is

$$A^* p = \frac{1}{2} Tr(H(\sigma\sigma^T p)) - \nabla \cdot (bp) \tag{2.7}$$

Where $H = (\partial_{x_i x_j})$ is the Hessian matrix, Tr evaluates the trace of a matrix, and $H(\sigma\sigma^T p)$ denotes the matrix multiplication of $H, \sigma, \sigma^T p$.

We take two simple example to show you how to get most probable phase portraits.

**Example 2.1:** Consider a scalar SDE with additive noise

$$dX_t = \alpha X_t dt + \beta dB_t \quad , \quad X(0) = x_0 \tag{2.8}$$

where $\alpha, \beta$ are real constants.

We solve the SDE first. We try for a solution of the form

$$X(t) = X_1(t)X_2(t) \tag{2.9}$$

where

$$\begin{cases} dX_1 = \alpha X_1 dX_1 \\ X_1(0) = 1 \end{cases} \tag{2.10}$$

and

$$\begin{cases} dX_2 = C(t)dt + D(t)dB_t \\ X_2(0) = x_0 \end{cases} \tag{2.11}$$

The function C(t), D(t) to be chosen. Then

$$\begin{aligned} dX &= X_2 dX_1 + X_1 dX_2 + 0 \cdot dt \\ &= \alpha X dt + C(t) X_1 dt + D(t) X_1 dB_t \\ &= \alpha X dt + \beta dB_t \end{aligned} \qquad (2.12)$$

For this, $C(t)=0$, $D(t)=\beta(X_1)^{-1}$ will work. Thus (2.11) reads

$$\begin{cases} dX_2 = \beta(X_1)^{-1} dB_t \\ X_2(0) = x_0 \end{cases} \qquad (2.13)$$

Now, we begin to solute $X_1$.

We set $Y_1 = \log X_1$, then using *Itô's fomula*, we have

$$\begin{aligned} dY_1 &= \frac{1}{X_1} dX_1 - \frac{1}{2X_1^2} \cdot 0 \cdot dt \\ &= \alpha dt \end{aligned} \qquad (2.14)$$

Thus

$$Y_1 = \alpha t + c, \ X_1 = C_0 e^{\alpha t} \qquad (2.15)$$

And we know that $X_1(0)=1$, it's easy to get $C_0=1$.

So we get

$$X_1 = e^{\alpha t} \qquad (2.16)$$

According to (2.13), we can get

$$X_2 = \beta \int_0^t e^{-\alpha s} dB_s + x_0 \qquad (2.17)$$

We conclude that

$$X(t) = X_1(t) X_2(t) = x_0 e^{\alpha t} + \beta \int_0^t e^{\alpha(t-s)} dB_s \qquad (2.18)$$

This is a Gaussian process. Once we have its mean and variance, we can then have its probability density function $p(x,t)$.

Its mean is

$$\mathbb{E} X(t) = e^{\alpha t} x_0 \qquad (2.19)$$

and the variance is

$$\begin{aligned}
Var(X(t)) &= \mathbb{E}\left[X(t) - e^{\alpha t}x_0\right]^2 \\
&= \mathbb{E}\left[\beta \int_0^t e^{\alpha(t-s)}dB_s\right]^2 \\
&= \beta \int_0^t e^{2\alpha(t-s)}ds = \frac{\beta}{2\alpha}\left[e^{2\alpha t} - 1\right]
\end{aligned} \qquad (2.20)$$

So the probability density function for the solution $X(t)$ is

$$\begin{aligned}
p(x,t) &= \frac{1}{\sqrt{2\pi}\sqrt{\frac{\beta}{2\alpha}}\sqrt{e^{2\alpha t}-1}} \exp\left(-\frac{(x-e^{\alpha t}x_0)^2}{\frac{\beta}{\alpha}(e^{2t}-1)}\right) \\
&= \frac{\sqrt{\alpha}}{\sqrt{\pi\beta}\sqrt{e^{2\alpha t}-1}} \exp\left(-\frac{\alpha(x-e^{\alpha t}x_0)^2}{\beta(e^{2\alpha t}-1)}\right)
\end{aligned} \qquad (2.21)$$

By setting $\partial_x p = 0$ or just observe $p(x,t)$, we obtain the maximizer at time t:

$$x_m(t) = e^{\alpha t}x_0 \qquad (2.22)$$

for every $x_0 \in \mathbb{R}^1$. Thus, the most probable dynamical system is

$$\dot{x}_m = \alpha x_m \qquad (2.23)$$

This is the same as the corresponding deterministic dynamical system $\dot{x} = \alpha x$, in the system with additive noise.

The phase portrait for the corresponding deterministic system and the most probable phase portrait are in Figure 2.1

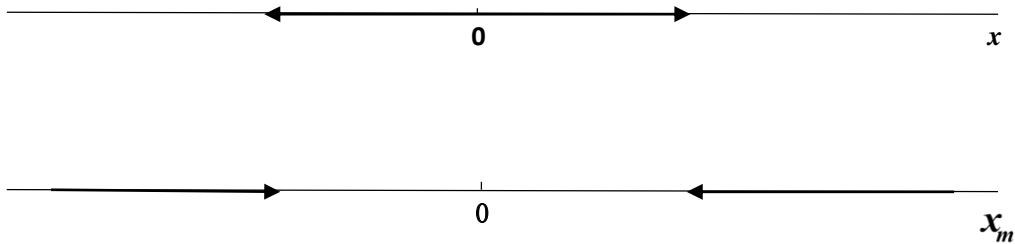

**Figure 2.1** Top: Phase portrait for the corresponding deterministic system (noise is absent) in Example 1.1. Bottom: Most probable phase portrait for Example 2.1.

**Example 2.2:** Consider a scalar SDE with multiplicative noise.

$$dX_t = \mu X_t dt + \sigma X_t dB_t \quad , \quad X(0) = x_0 \qquad (2.24)$$

where $\mu, \sigma$ are real constants.

We solve the SDE first. Let

$$Y_t = \log X_t \tag{2.25}$$

and so

$$dY_t = d(\log X_t) = \frac{dX_t}{X_t} - \frac{1}{2}\frac{\sigma^2 X_t^2 dt}{X_t^2} \quad \text{by Itô's fomula}$$
$$= (\mu - \frac{\sigma^2}{2})dt + \sigma dB_t \tag{2.26}$$

Consequently

$$X_t = x_0 e^{(\mu - \frac{\sigma^2}{2})t + \sigma B_t} \tag{2.27}$$

Then let us find out its probability distribution function $F(x,t)$ and then its probability density function $p(x,t)$.

The state 0 is an equilibrium state. When the initial state $x_0$ is positive (negative), the solution is also positive (negative). The distribution function is calculated as

$$\begin{aligned}F(x,t) &= P(X(t) \leq x) = P(x_0 e^{(\mu - \frac{\sigma^2}{2})t + \sigma B_t} \leq x) \\ &= P(e^{\sigma B_t} \leq \frac{x}{x_0} e^{-(\mu - \frac{\sigma^2}{2})t}) \\ &= P(B_t \leq \frac{1}{\sigma}(\log(\frac{x}{x_0}) - (\mu - \frac{\sigma^2}{2})t)) \\ &= \int_0^{\frac{1}{\sigma}(\log(\frac{x}{x_0}) - (\mu - \frac{\sigma^2}{2})t)} \frac{1}{\sqrt{2\pi t}} e^{-\frac{\xi^2}{2t}} d\xi \end{aligned} \tag{2.28}$$

For $x_0 \neq 0$, the distribution function F is non-zero only if x and $x_0$ have the same sign (i.e., $\frac{x}{x_0}$ is positive). The probability density function is then

$$p(x,t) = \partial_x F(x,t) = \frac{1}{\sigma x}\frac{1}{\sqrt{2\pi t}}\exp(-\frac{(\log(\frac{x}{x_0}) - (\mu - \frac{\sigma^2}{2})t)^2}{2\sigma^2 t}) \tag{2.29}$$

Setting $\partial_x p = 0$, we obtained

$$x_m(t) = x_0 e^{(\mu - \frac{3}{2}\sigma^2)t} . \tag{2.30}$$

It can be checked that this is indeed the maximizer of $p(x,t)$ at time t. Thus, the most probable dynamical system is

$$\dot{x}_m = (\mu - \frac{3}{2}\sigma^2)x_m \tag{2.31}$$

## 3. Motivation

As we see in the first section, the methods to get the most probable phase portrait can be comment as the following four steps.

**Step1:** figure out the solution of the stochastic differential equations;

**Step2:** figure out the probability density function $p(x,t)$ of the solution;

**Step3:** figure out the most probable solution by setting $\partial_x p = 0$;

**Step4:** plot the most probable phase portrait.

For many simple stochastic differential equations, we can get the most probable phase portraits easily by these four steps.

However, for some more complicated SDEs, solving the analytical solution is not easy at all and thus we cannot get the MPPP[1] easily. Therefore, it is significant to find numerical ways to get the MPPP. This passage will show you our numerical method later.

## 4. Numerical method to get most probable phase portraits

According to the four steps listed in the former section, we will follow the three steps to get the most probable phase portrait by MATLAB simulation:

**Step1:** plot many discretized paths of the solution to the stochastic differential equations;
**Step2:** get the most probable points at every discrete t values.
**Step3:** plot the discretized most probable phase portraits.
Then we will follow the three steps to introduce you our method. And we will take the SDE in Example 1.1 as example.

### 4.1 Numerical simulation of SDEs -- Euler-Maruyama method
A scalar, autonomous SDE can be written in integral form as

$$X_t = x_0 + \int_0^t f(X_s,s)ds + g(X_s,s)dB_s \tag{4.1}$$

---
[1] MPPP: abbreviation of Most Probable Phase Portraits.

```matlab
%% EM Euler-Maruyama method on linear SDE
%
% SDE is dX = alpha*X*dt + beta*dB , X(0) = Xzero,
% where alpha =1, beta=1 and Xzero = 1.
%
% Discretized Brownian path over [0,1] has dt = 2^(-7).
% Euler-Maruyama uses timestep dt.
%% I: initial datas
T = 1; % right point of the time interval
Xzero = 1; % initial value of x
M = 2^15; % number of paths
N = 2^7*T; % number of steps in every path
randn('state',M*N)
dt = T/N; % discretized time step
dB = sqrt(dt)*randn(M,N); % Brownian increments
X = ones(M,N); % preallocate for efficiency
%% II: figure out the value of X in every discretized time of every path
for i=1:M
   Xtemp = Xzero;
   for j = 1:N
       Xtemp = Xtemp + Xtemp*dt + dB(i,j);
       X(i,j) = Xtemp;
   end
end
%% III: plot the paths
for i=1:M
   plot([0:T/N:T],[Xzero,X(i,:)]) ,hold on
end
xlabel('t','FontSize',12)
ylabel('X','FontSize',16,'Rotation',0,'HorizontalAlignment','right'
)
```

**Listing 1**:   *M-file* Xpaths.m

where *f* and *g* are scalar functions and the initial condition is $X(0) = x_0$.

It is usual to rewrite (4.1) in differential equation form as

$$dX_t = f(X_t,t)dt + g(X_t,t)dB_t \quad , \quad X(0) = x_0 \qquad (4.2)$$

To apply a numerical method to (4.2) over [0, T], we first discretize the interval. Let $\Delta t = \dfrac{T}{N}$ for some positive integer N, and the partition is

$$0 = t_0 < t_1 < \cdots < t_N = T$$

Our numerical approximation to $X(\tau_j)$ will be denoted $X_j$. The Euler–Maruyama (EM) method takes the form

$$X_j = X_{j-1} + f(X_{j-1}, t_{j-1})\Delta t + g(X_{j-1}, t_{j-1})(B(t_j) - B(t_{j-1})), \quad j = 1, 2, \ldots, N \quad (4.3)$$

In this article, we will use discretized Brownian paths in [2, 527] to generate the increments needed $B(t_j) - B(t_{j-1})$ in (4.3).

We will apply the EM method to the linear SDE in Example 1.2

$$dX_t = \alpha X_t dt + \beta dB_t, \quad X(0) = x_0 \quad (4.4)$$

and we set $\alpha = \beta = 1$, $x_0 = 1$. You can see the code in Listing 1.

Run the code, we can get $2^{15}$ discretized paths of X, just as Figure 4.1 shows.

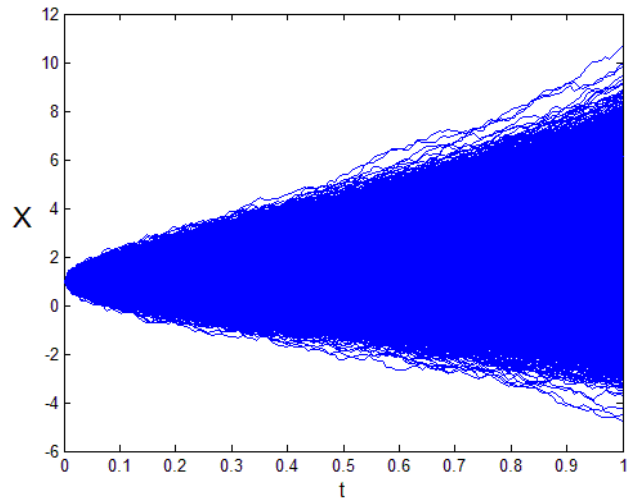

**Figure 4.1:** EM approximation of the solution for (4.4), from Xpaths.m

### 4.2 Numerical simulation of MPPP

Now we have get the many simulation paths of SDE. The next step is to figure out the most probable points at every discrete t values.

We use the function ksdensity in MATLAB to achieve our goal. For every discretized time, there are M points in M paths of X and we use function ksdensity in MATLAB to get the kernel smoothing density estimation of X. Then we can easily find the most probable points in every time, which is corresponding to the highest points of kernel smoothing density estimation at every time.

From Example 1.1, we can see that the true MPPP is $x_m(t) = e^t x_0$. And we plot the true MPPP in blue line to see if the simulation result matches well.

```matlab
%% Ploting MPPP
%
% SDE is dX = alpha*X*dt + beta*dB , X(0) = Xzero,
% where alpha =1, beta=1 and Xzero = 1.
%
% Discretized Brownian path over [0,1] has dt = 2^(-7).
% Euler-Maruyama uses timestep dt.
%% I: initial data
T = 1; % right point of the time interval
Xzero = 1; % initial value of x
M = 2^15; % number of paths
N = 2^8*T; % number of steps in every path
randn('state',M*N)
dt = T/N; % discretized time step
dB = sqrt(dt)*randn(M,N); % Brownian increments
X = ones(M,N); % preallocate for efficiency
%% II: figure out the value of X in every discretized time of every path
for i=1:M
   Xtemp = Xzero;
   for j = 1:N
       Xtemp = Xtemp + Xtemp*dt + dB(i,j);
       X(i,j) = Xtemp;
   end
end
%% III: figure out the most probable points at every discrete t values.
MPPP = zeros(1,N); % preallocate for efficiency
for t=1:N
   Xemt = X(:,t);
   [h,xi] = ksdensity(Xemt);
   [hmax,I] = max(h);
   MPPP(t) = xi(I);
end
%% IV: plot the real most probable phase portraits
tt = [0:T/N:T];
MPPP_true = exp(tt);
plot(tt,MPPP_true,'b-'),hold on
MPPP_err = abs(MPPP(N)-MPPP_true(N))/MPPP_true(N)
%% V: plot the discretized most probable phase portraits.
plot([0:T/N:T],[Xzero,MPPP],'r-')
xlabel('t','FontSize',12)
ylabel('MPPP','FontSize',16)
```

**Listing 2**:　　*M-file* MPPP.m

Running the code in Listing 2, we can get **Figure 4.2**.

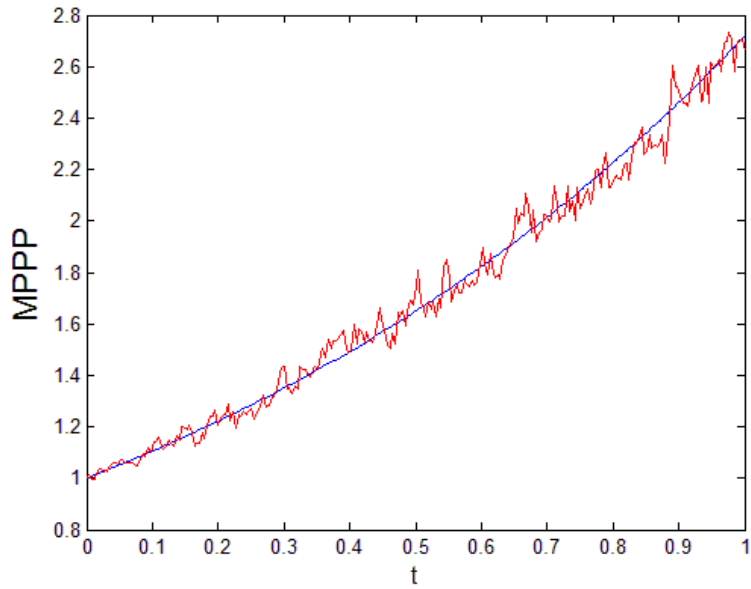

**Figure 4.2:** True MPPP(blue line) and simulated MPPP(red line) for (4.4), from MPPP.m

When we set T to be 2, 4, 6, 8, we can get Figure 4.3

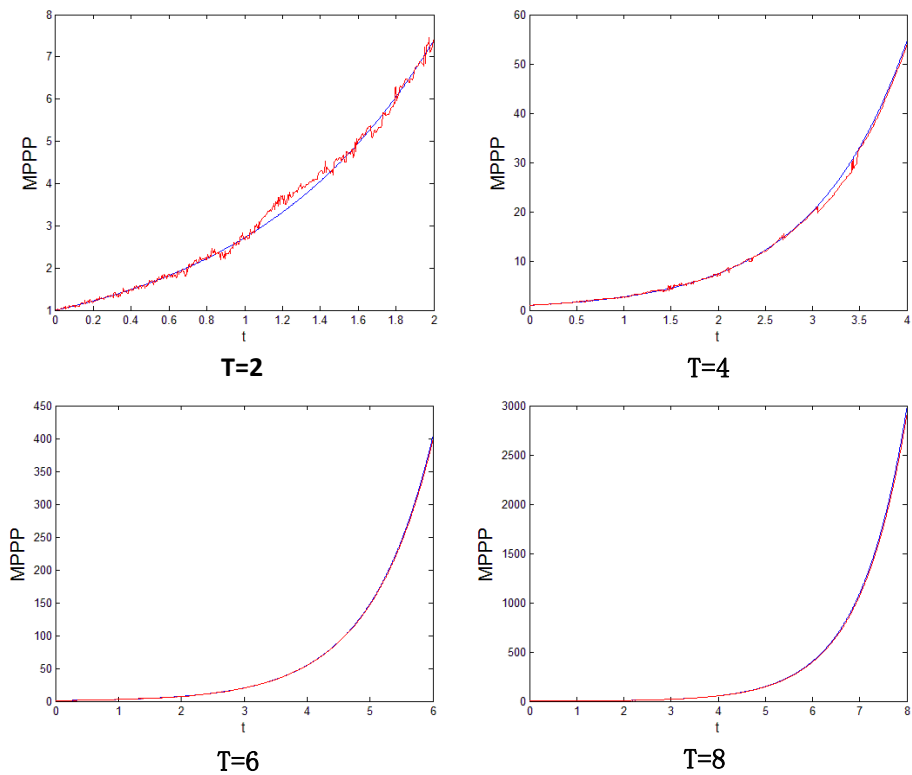

**Figure 4.3:** MPPPs when T = 2, 4, 6, 8 respectively

```
%% III: figure out the value of X in every discretized time of
every path
for i=1:M
   Xtemp = Xzero;
   for j = 1:N
      Xtemp = Xtemp + f(Xtemp,t)*dt + g(Xtemp,t)dB(i,j);
      X(i,j) = Xtemp;
   end
end
```

**Listing 3**: *M-file* MPPP.m

From Figure 4.3, we can see that our method works very well.

Furthermore, we have to note that the numerical method is broadly applicable for all SDEs.

**For one-dimension SDEs**

$$dX_t = f(X_t,t)dt + g(X_t,t)dB_t \quad , \quad X(0) = x_0 \tag{4.5}$$

We just need to correspondingly change the second part of the MATLAB code in Listing 2 to be Listing 3.

**For two-dimension SDEs**, this method could also work, take (4.6) for example.

$$\begin{cases} dX_t = -Y_t dt + B_t^1 \\ dY_t = X_t dt + B_t^2 \end{cases} \tag{4.6}$$

We can get its numerical simulation result by the code in Listing 4.

Running the code, we can get the MPPP of X and Y, see in Figure 4.4.

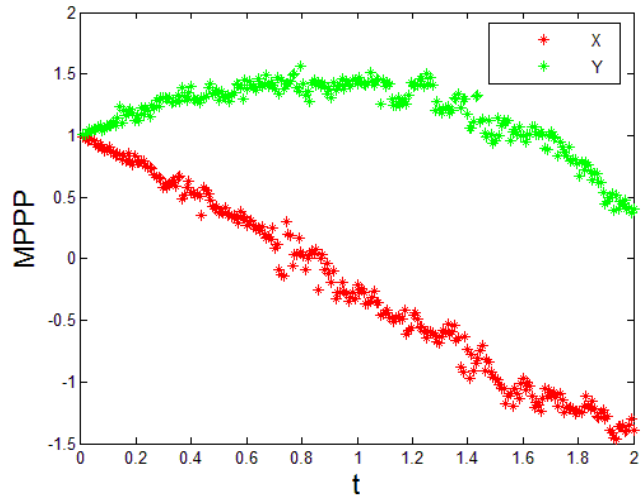

**Figure 4.4:** simulated MPPP of X and Y for (4.6), from MPPP_2.m

```matlab
%% Ploting MPPP of 2-dimension SDE
% SDE is dX = Y*dt + dB1
%        dY = X*dt + dB2 , X(0)=Xzero, Y(0)=Yzero
% Discretized Brownian path over [0,2] has dt = 2^(-7).
% Euler-Maruyama uses timestep dt.
%% I: initial datas
T = 2; % right point of the time interval
Xzero = 1; Yzero = 1; % initial value of x
M = 2^15; % number of paths
N = 2^7*T; % number of steps in every path
randn('state',2*M*N)
dt = T/N; % discretized time step
dB1 = sqrt(dt)*randn(M,N); dB2 = sqrt(dt)*randn(M,N);% Brownian increments
X = ones(M,N); Y = ones(M,N); % preallocate for efficiency
%% II: figure out the value of X in every discretized time of every path
for i=1:M
   Xtemp = Xzero; Ytemp = Yzero;
   for j = 1:N
      Xtemp = Xtemp - Ytemp*dt - dB1(i,j);
      Ytemp = Ytemp + Xtemp*dt + dB2(i,j);
      X(i,j) = Xtemp; Y(i,j) = Ytemp;
   end
end
%% III: figure out the most probable points at every discrete t values.
MPPPx = zeros(1,N) ; MPPPy = zeros(1,N); % preallocate for efficiency
for t=1:N
   Xemt = X(:,t); Yemt = Y(:,t);
   [hx,xi] = ksdensity(Xemt); [hy,yi] = ksdensity(Yemt);
   [hxmax,Ix] = max(hx); [hymax,Iy] = max(hy);
   MPPPx(t) = xi(Ix); MPPPy(t) = yi(Iy);
end
%% IV: plot the discretized most probable phase portraits.
plot([0:T/N:T],[Xzero,MPPPx],'r*',[0:T/N:T],[Yzero,MPPPy],'g*')
xlabel('t','FontSize',16)
ylabel('MPPP','FontSize',16)
%plot([Xzero,MPPPx],[Yzero,MPPPy],'r*')
%xlabel('X','FontSize',16)
%ylabel('Y','FontSize',16,'Rotation',0,'HorizontalAlignment','right')
```

**Listing 4**: *M-file* MPPP_2.m

And the phase portraits in X-Y phase plain is

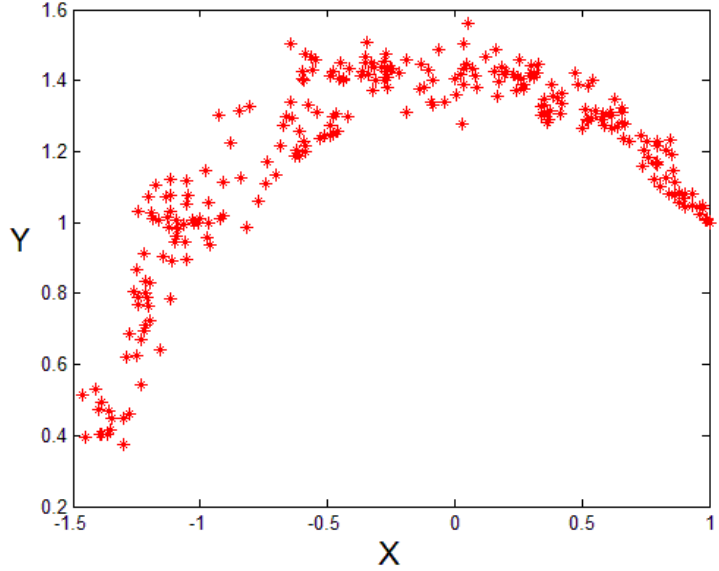

**Figure 4.5:** simulated MPPP in X-Y plain for (4.6), from MPPP_2.m

Just as we noted in one-dimension SDEs, the method is also applicable for the general cases

$$\begin{cases} dX_t = f_1(X_t, Y_t, t)dt + g_1(X_t, Y_t, t)dB_t^1 \\ dY_t = f_2(X_t, Y_t, t)dt + g_2(X_t, Y_t, t)dB_t^2 \end{cases} \quad (4.7)$$

**For three-dimension and more higher dimension SDEs**, obviously the method can also work.

## 5. Summary

Recalling all of the above, we have established a better way to get phase portraits of SDEs. For simple SDEs, it's easy to get its analytical solution by the four steps in section 3. For more complicated SDEs, we can use the numerical simulation method in section 4 and the numerical method is broadly applicable in all SDEs.

Though it seems that we have solved the question to get MPPP, there are still many problems to be solved later. The method above is broadly applicable for all SDEs, but for some SDEs, it cannot fit the real result very well or the MPPP is very unstable because of the stochastic factor. Thus it still have spaces for further improvement. Moreover, analyses the convergence of the methods is another problem to be solved, which may also be helpful for further improvement of the method.